\documentclass[12pt,reqno,a4paper]{amsart}
\usepackage{amsaddr}
\usepackage{amsthm,amsmath,amssymb,eucal}
\usepackage{geometry}
\usepackage{resizegather}
\usepackage{graphicx}
\usepackage{xcolor}
\usepackage[normalem]{ulem}
\usepackage[latin1]{inputenc}
\usepackage{enumerate}
\usepackage{booktabs}

\newcommand{\ee}{\mathrm{e}}
\newcommand{\D}{\mathrm{d}}
\newcommand{\C}{\mathbb{C}}
\newcommand{\N}{\mathbb{N}}

\newcommand{\Z}{\mathbb{Z}}
\newcommand{\DD}{\mathcal{D}}

\newcommand{\PP}{\mathcal{P}}

\newcommand{\Heff}{H_{\mathrm{eff}}}

\newtheorem{claim}{Claim}[section]
\newtheorem{theorem}[claim]{Theorem}
\newtheorem{proposition}[claim]{Proposition}

\newtheorem{corollary}[claim]{Corollary}
\newtheorem{remark}[claim]{Remark}
\newtheorem{remarks}[claim]{Remarks}

\begin{document}

\title[magnetic lattice graph with a circulant vertex coupling]{A new butterfly in the meadow: magnetic lattice graphs with a circulant vertex coupling}

\author{Pavel Exner}

\address{{Doppler Institute for Mathematical Physics and Applied Mathematics, Czech Technical University, B\v rehov{\'a} 7, 11519 Prague}
{\rm and} {Department of Theoretical Physics, Nuclear Physics Institute, Czech Academy of Sciences, 25068 \v{R}e\v{z} near Prague, Czechia}
}

\email{exner@ujf.cas.cz}

\author{Jan Peka\v{r}}

\address{{Faculty of Mathematics and Physics, Charles University, V Hole\v{s}ovi\v{c}k\'ach 2, 18040 Prague}
{\rm and} {Department of Theoretical Physics, Nuclear Physics Institute, Czech Academy of Sciences, 25068 \v{R}e\v{z} near Prague, Czechia}
}

\email{honzapekar28@gmail.com}

\date{\today}

\maketitle

\begin{flushright}
\textit{To Uzy Smilansky, one of quantum graph prophets, \\ on the occasion of his 85th birthday}
\end{flushright}

\begin{abstract}
We analyze the spectrum of quantum square lattice graph exposed to a homogeneous magnetic field with a rational flux through plaquette, assuming all nontrivial cases of the vertex coupling which is circulant and scale-invariant, with a particular attention to those of them that give rise to a dynamics violating the time reversal invariance. Except for two cases with partially decoupled edges, the remaining twelve split into two spectral classes; the second one offers a hitherto unknown butterfly pattern.
\end{abstract}

\section{Introduction} 
\setcounter{equation}{0}

The idea of quantum mechanics with a graph as configuration space was proposed by Linus Pauling ninety years ago \cite{Pa36} but then happily forgotten for half a century. Its second birth came in the eighties, motivated primarily, but not exclusively, by the use of such operators as models of quantum nanostructures which then started entering the scene. It took time and vision, however, to reveal the true potential of the concept, and the birthday boy was one of the main movers in the field; to limit ourselves to a small sample one can mention, e.g., his work of quantum chaos on graphs \cite{KS97}, trace formul{\ae} \cite{KS99}, inverse problems \cite{GS01}, irreversible dynamics \cite{Sm04}, nodal domains on graphs \cite{BBS12}, and many others.

While it is generally known that the self-adjointness requirement allows to match wave functions in the graph vertices by many different ways, most results in the literature concern the simplest possibility, usually called Kirchhoff, or at most the one-parameter family of the $\delta$-couplings~\cite[Sec.~1.4.4]{BK13}. To look into other couplings, one needs a motivation. A possible one comes from the attempt to use quantum graphs to model the anomalous Hall effect \cite{SK15, SV23}. In order to get the result, one had there to impose a preferred direction of motion on the graphs edges, a step that cannot be justified from the first principles. One can, however, do the same by choosing a vertex coupling which violates the time-reversal invariance.

One of the examples we treated concerned a square lattice with a time asymmetric vertex coupling exposed to a homogeneous magnetic field \cite{BEL23}. It is of a particular interest because the vertex coupling and the magnetic field represent two competing mechanisms to violate the invariance: what was found was that the former dominates at low energies and the latter in the high-energy asymptotics exhibiting the well-known butterfly pattern. Recall that interesting spectral properties of systems having two characteristic lengths, possibly incommensurate, was foretold already by Azbel \cite{Az64} but it was only the Hofstadter illustrative result \cite{Ho76} which managed to ignite the interest to this problem.

The fact that the Hofstadter butterfly dominates in the model of \cite{BEL23} at high energies is easy to explain: the time-reversal invariance violation is there energy dependent and it vanishes asymptotically. This naturally inspires the question how such a system would behave if the coupling would be \emph{scale-invariant} which is the problem we address in this paper. It is well known \cite{KS03} that the most general self-adjoint vertex coupling in a quantum graph exposed to a magnetic field (and free of other interactions) is given by the matching condition
\begin{equation} \label{eq:coupling}
(U - I)\,\psi(v) + i\,(U + I)\,(\DD\psi)(v) = 0
\end{equation}
for the boundary values of the wave function $\psi$ at the vertex $v$, where $\DD=\frac{\D}{\D x}-iA_e$ with $x$ being the coordinate parametrizing an edge $e$ rooted in the vertex and $A_e$ the component of the vector potential tangent to $e$, and $U$ is an $n\times n$ unitary matrix with $n=\mathrm{deg}\,v$. The scale invariance means the absence of the Robin component of the coupling in the terminology of ~\cite[Thm.~1.4.4]{BK13}, so that the matrix $U$ is at the same time Hermitean. To make the problem tractable, we also assume that the coupling is \emph{circulant}, i.e. invariant with respect to a cyclic permutation of the edges which in case of a planar graph makes a clear sense.

\begin{figure}[htb!]
		\centering
		\includegraphics[scale=0.4]{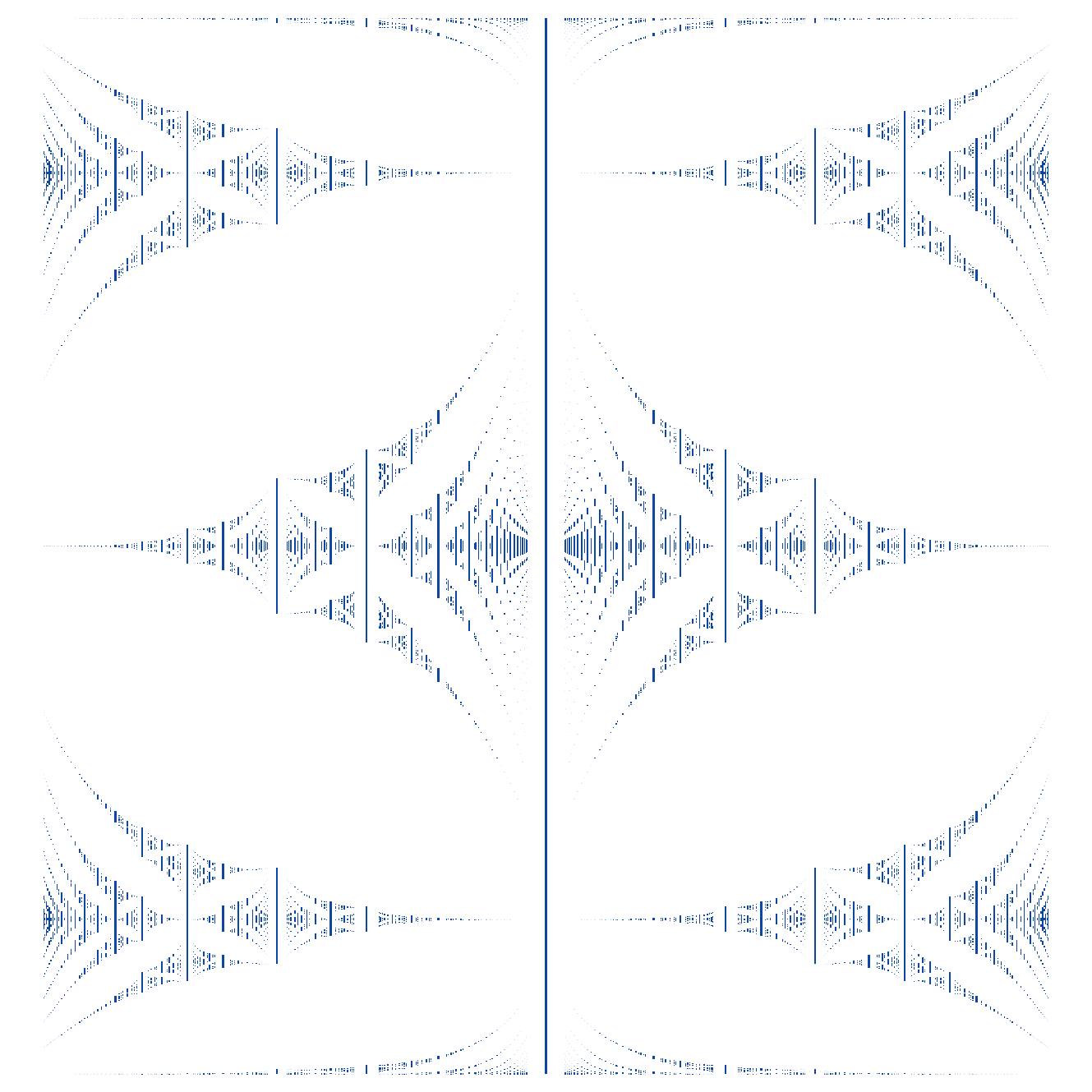}
		\caption{Spectrum of a lattice with class II coupling, cf. Theorem~\ref{thm:class}.}
	\label{newbutterfly}
	\end{figure}

We are going to list all the couplings with these properties, altogether fourteen cases, and show that, apart from the simple situation when the lattice edges are pairwise decoupled, the spectrum for magnetic fields, rational in terms of the flux per plaquette, splits into two classes. One is the classical Hofstadter butterfly, the other is a pattern which have not yet been observed, cf.~Figure~\ref{newbutterfly}.

\section{Preliminaries: a class of circulant matrices} \label{s:matrices}
\setcounter{equation}{0}

In accordance with the motivation mentioned in the introduction, we are interested in a particular class of circulant matrices which are simultaneously Hermitean and unitary, which means that such a $U$ must be idempotent, $U^2 = I$. Let us recall some of their properties.

\subsection{The general case}

An $n \times n$ matrix $U$ is called \emph{circulant} if it is of the form
 $$ 
		U = \begin{pmatrix} c_1 & c_2 & c_3 & {\dots} & c_n \\
			c_n & c_1 & c_2 & {\dots} & c_{n-1} \\
			{\vdots} & c_n & c_1 & {\ddots} & {\vdots} \\
			c_3 & {\dots} & {\ddots} & {\ddots} & c_2 \\
			c_2 & c_3 & {\dots} & c_n & c_1 \\
		\end{pmatrix},
 $$ 
being uniquely determined by its generating vector $(c_1, c_2, c_3,{\dots}, c_n)^{\top}$, $c_j \in \mathbb{C}$. Such matrices are diagonalized by the discrete Fourier transform, $D = FCF^*$, where $F_{rs}=\frac{1}{\sqrt{n}}\, \omega^{-(r-1)(s-1)}$, $\,r,s=1,\dots,n$ and ${\omega} := \ee^{2{\pi}i/n}$ is the $n$th complex root of unity. All $n \times n$ circulant matrices have the same eigenvectors,
    \begin{equation} \label{eq:efs}
		\chi_m = \frac{1}{\sqrt{n}}\{1,\omega^m,\omega^{2m},\dots,\omega^{(n-1)m}\}^\mathsf{T}, \quad m=0,1,2,\dots,n-1.
	\end{equation}
The eigenvalues of such matrices are expressed through the generating vector, $\lambda_j = \sum_{a=1}^{n} c_a{\omega}^{j(a-1)}$, and vice versa, $c_j = \frac{1}{n}\sum_{a=1}^{n} \lambda_a{\omega}^{-a(j-1)}$, where the indices are understood modulo $n$, i.e. $\lambda_n=\lambda_0$, etc. The condition $U^2 = I$ implies that the eigenvalues $\lambda_j$ of $U$ can be only $\pm 1$, and we will write them as an ordered set $D := (\lambda_0, \lambda_1, \lambda_2,{\dots}, \lambda_{n-1})$; we also introduce the sets $S_\pm=\{m:\, \lambda_m=\pm 1\}$.

The relations between the generating vector and the eigenvalues show that each circulant matrix of this class in uniquely determined by the set $D$, or equivalently by $S:= S_+$. This means that there are $2^n$ of them of which we will exclude the trivial situations where $U=\mp I$ for which the edges in \eqref{eq:coupling} are fully decoupled with Dirichlet and (magnetic) Neumann endpoints, respectively. What is important is that the order of $D$ matters; a reordering yields a different circulant matrix. For instance, a cyclic permutation of $D$ given by $\lambda'_k=\lambda_{k+r}$ (mod $n$) is easily seen to produce matrix $U'$ the entries of which are $U'_{jk}= \omega^{r(k-j)} U_{jk}$.

Let us further recall that according to \cite{ET21} a quantum graph with vertex coupling given by \eqref{eq:coupling} is \emph{time-reversal invariant} if and only if the matrix $U$ is transpose invariant, which is further equivalent to the invariance with respect to the complex conjugation; for circulant matrices this means $c_j = c_{n-j+2}$, $j=2,\dots,\lfloor\frac{n+1}{2}\rfloor$.

\subsection{The $n=4$ case} \label{ss:nfour}

Let us now specify the above general claims to the $4\times 4$ matrices which are of our interest here. It follows from the above discussion that there are fourteen nontrivial possibilities looking explicitly as follows:
\begin{align}
& U_{\{0\}} = \frac{1}{2}{\scriptsize \begin{pmatrix} -1& 1 & 1 &1 \\
			                            1 & -1 & 1 & 1 \\
                            			1 & 1 & -1 & 1 \\
                            			1 & 1 & 1 & -1 \\
                        \end{pmatrix}}, \quad
U_{\{1,2,3\}} = \frac{1}{2}{\scriptsize \begin{pmatrix} 1& -1 & -1 &-1 \\
			                            -1 & 1 & -1 & -1 \\
                            			-1 & -1 & 1 & -1 \\
                            			-1 & -1 & -1 & 1 \\
                        \end{pmatrix}} \label{Kirch,antiKirch} \\[.3em]
 & U_{\{1\}} = \frac{1}{2}{\scriptsize \begin{pmatrix} -1& -i & -1 &i \\
			                            i & -1 & -i & -1 \\
                            			-1 & i& -1 & -i \\
                            			-i & -1 & i & -1 \\
                        \end{pmatrix}}, \quad
U_{\{0,2,3\}} = \frac{1}{2}{\scriptsize \begin{pmatrix} 1& i & 1 &-i \\
			                            -i & 1 & i & 1 \\
                            			1 & -i& 1 & i \\
                            			i & 1 & -i & 1 \\
                        \end{pmatrix}}  \label{single1}  \\[.3em]
& U_{\{2\}} = \frac{1}{2}{\scriptsize \begin{pmatrix} -1& -1 & 1 &-1 \\
			                            -1 & -1 & -1 & 1 \\
                            			1 & -1 & -1 & -1 \\
                            			-1 & 1 & -1 & -1 \\
                        \end{pmatrix}}, \quad
U_{\{0,1,3\}} = \frac{1}{2}{\scriptsize \begin{pmatrix} 1& 1 & -1 &1 \\
			                            1 & 1 & 1 & -1 \\
                            			-1 & 1 & 1 & 1 \\
                            			1 & -1 & 1 & 1 \\
                        \end{pmatrix}}  \label{single2} \\[.3em]
& U_{\{3\}} = \frac{1}{2}{\scriptsize \begin{pmatrix} -1& i & -1 &-i \\
			                            -i & -1 & i & -1 \\
                            			-1 & -i& -1 & i \\
                            			i & -1 & -i & -1 \\
                        \end{pmatrix}}, \quad
U_{\{0,1,2\}} = \frac{1}{2}{\scriptsize \begin{pmatrix} 1& -i & 1 &i \\
			                            i & 1 & -i & 1 \\
                            			1 & i& 1 & -i \\
                            			-i & 1 & i & 1 \\
                        \end{pmatrix}}  \label{single3} \\[.3em]
& \hspace{-1.5em} U_{\{0,1\}} = \frac{1}{2}{\scriptsize \begin{pmatrix} 0 & 1-i & 0 &1+i \\
			                            1+i & 0 & 1-i & 0 \\
                            			0 & 1+i & 0 & 1-i \\
                            			1-i & 0 & 1+i & 0 \\
                        \end{pmatrix}}, \,
U_{\{2,3\}} = \frac{1}{2}{\scriptsize \begin{pmatrix} 0 & -1+i & 0 &-1-i \\
			                            -1-i & 0 & -1+i & 0 \\
                            			0 & -1-i & 0 & -1+i \\
                            			-1+i & 0 & -1-i & 0 \\
                        \end{pmatrix}}  \label{double_side} \\[.3em]
& \hspace{-1.5em} U_{\{0,3\}} = \frac{1}{2}{\scriptsize \begin{pmatrix} 0 & 1+i & 0 &1-i \\
			                            1-i & 0 & 1+i & 0 \\
                            			0 & 1-i & 0 & 1+i \\
                            			1+i & 0 & 1-i & 0 \\
                        \end{pmatrix}}, \,
U_{\{1,2\}} = \frac{1}{2}{\scriptsize \begin{pmatrix} 0 & -1-i & 0 &-1+i \\
			                            -1+i & 0 & -1-i & 0 \\
                            			0 & -1+i & 0 & -1-i \\
                            			-1-i & 0 & -1+i & 0 \\
                        \end{pmatrix}}  \label{double_middle} \\[.3em]
& \qquad\qquad U_{\{0,2\}} = {\small \begin{pmatrix} 0 & 0 & 1 &0 \\
			                    0 & 0 & 0 & 1 \\
                            	1 & 0 & 0 & 0 \\
                            	0 & 1 & 0 & 0 \\
                \end{pmatrix}}, \quad
U_{\{1,3\}} = {\small \begin{pmatrix} 0 & 0 & -1 &0 \\
			                    0 & 0 & 0 & -1 \\
                            	-1 & 0 & 0 & 0 \\
                            	0 & -1 & 0 & 0 \\
                \end{pmatrix}} \label{pair_decoupl}
\end{align}
As indicated, all of them have a common set of eigenvectors; relation \eqref{eq:efs} simplifies in this case to
 $$ 
		\chi_m = \frac12\{1,i^m,i^{2m},i^{3m}\}^\mathsf{T}, \quad m=0,1,2,3\,;
 $$ 
the spectral decomposition of a given $U_S$ can be then written explicitly as
 $$ 
		U_S= \sum_{m=0}^3 \lambda_m (\chi_m,\cdot)\chi_m = E_+ - E_-,
 $$ 
where the $E_\pm$ are projections onto $\mathrm{span}\{\chi_m:\, \lambda_m=\pm 1\}$. In each pair \eqref{Kirch,antiKirch}--\eqref{pair_decoupl} the matrices differ by sign, and moreover, the pairs \eqref{single1} and \eqref{single3}, as well as \eqref{double_side} and \eqref{double_middle}, are complex conjugated.
\begin{remark} \label{maxviol}
{\rm Out of the fourteen matrices \eqref{Kirch,antiKirch}--\eqref{pair_decoupl} eight are not invariant with respect to transposition, hence the corresponding vertex couplings give rise to the dynamics violating the time-reversal invariance. One may ask what is the measure of this violation. A natural characterization is given by the norm of the difference $U_S-U_S^\mathsf{T}$. It is easy to check that each of these matrices has eigenvalues zero of multiplicity two and $\pm 2$, meaning that $\|U_S-U_S^\mathsf{T}\|=2$ so that the violation is maximal saturating the unitarity bound. In contrast to the general circulant vertex couplings considered in \cite{ET21}, the scale-invariant ones can be thus with respect to the time reversal extremal only, either such a coupling preserves the invariance or it violates it maximally.}
\end{remark}

\section{Statement of the problem} \label{s:statement}
\setcounter{equation}{0}

\subsection{The model} \label{ss:model}

We consider an infinite square lattice graph; without loss of generality we may suppose that its vertices are the points of $\Z^2$. The graph is the configuration space of a charged particle exposed to a homogeneous magnetic field of strength $B$ perpendicular to the lattice plane. The Hamiltonian thus acts as $\psi_j \mapsto \DD^2\psi_j$ on the $j$th edge, where $\DD=\frac{\D}{\D x}-iA$ is the \emph{outward} magnetic quasi-derivative; we denote it as $-\Delta_{B,U}$ with $U$ indicating the vertex coupling. We work in the Landau gauge choosing the vector potential as $\vec A=(0,Bx,0)$, and use the rational units, $\hbar=2m=e=c=1$, in which the magnetic flux quantum is $\Phi_0=2\pi$; we assume that the flux per plaquette is a rational multiple of this quantity, $\Phi= 2\pi\frac{p}{q}$, with coprime $q\ge 2$ and $p\in\{1,\dots,q-1\}$. In other words, the field strength is $B=2\pi\frac{p}{q}$, vertical edges in the $m$th column carry the tangential potential equal to $Bm$ while the horizontal ones carry none, as in~\cite{BEL23}.

At the lattice nodes the wave functions satisfy the conditions \eqref{eq:coupling} in which $U=U_S$ is one of the matrices \eqref{Kirch,antiKirch}--\eqref{pair_decoupl}. This is the crucial assumption; any vertex coupling of this type is scale-invariant which means that the matching condition decouple and can be written in the projection form,
\begin{equation} \label{eq:projform}
E_-\,\psi(v) = 0, \quad E_+\,(\DD\psi)(v) = 0
\end{equation}
with the boundary values laying in $E_+\C^4$ and the magnetic quasi-derivatives values in $E_-\C^4$. Other consequences of the scale invariance are that the on-shell vertex scattering matrix \cite[Sec.~2.1.1]{BK13} is energy-independent, $S(k) = U$ for all $k$, and that the negative spectrum of any $-\Delta_{B,U_S}$ is empty by \cite[Thm.~2.6]{BET22}. Our aim is to find the positive spectrum for all the indicated vertex couplings.

\subsection{Duality} \label{ss:duality}

The scale-invariant character of the coupling allows us, in contrast to \cite{BEL23}, to find an equivalent formulation of the problem in terms of an effective difference operator. This is a tool often used in such situations; we will comment on various related problems in Sec.~\ref{ss:related} below.

Let us denote by $d$ the dimension of the positive spectral subspace of the matrix determining the coupling, $d:=|S|=\dim E_+ \in\{1,2,3\}$. Given the splitting of the boundary values and quasi-derivative values at a vertex $v$ given by \eqref{eq:projform}, the former can be written as $\psi(v)= \sum_{a\in S} \phi_a(v)\chi_a$ with complex coefficients $\phi_a(v)$. For further purposes we note that the entries of the eigenvectors $\chi_m$ differ in the `left-right' and `up-down' pairing, $\tilde\jmath := j+2 \pmod 4$, at most by sign
\begin{equation} \label{eq:reversal}
(\chi_m)_{\tilde\jmath} = i^{2m} (\chi_m)_j = (-1)^m (\chi_m)_j .
\end{equation}
Next we define for $j=1,\dots,4$ and $a,b\in S$ the $d\times d$ \emph{hopping matrices} with entries
\begin{equation} \label{eq:Tj}
\bigl(T_j\bigr)_{ab} := (-1)^{b}\; i^{\,(b-a)(j-1)} .
\end{equation}
Each matrix $T_j$ is rank one: to see it, take vectors $v_j\in\C^d$ with components $(v_j)_a=i^{a(j-1)}$, then one can write the action of $T_j$ as $(T_j\phi)_a = (\bar v_j)_a \sum_{b\in S} (-1)^b (v_j)_b\phi_b$. Furthermore, it is straightforward to check that we have  $T_3 = T_1^*$ and $T_4 = T_2^*$.
\begin{remarks} \label{hop_differences}
{\rm (a) The hopping matrices depend on the choice of $U$ through the indices contained in $S$ which enter in two different ways, through the differences \mbox{$b-a\:(\bmod\;4)$} and through the parities $(-1)^a$. For $d=2$ those are interconnected leading to two different direction types; note that the diagonal of $T_j$ is $j$-independent and the direction enters only through the off-diagonal factor $i^{(b-a)(j-1)}$. For $b - a = 2$ this equals $(-1)^{j-1}$, which depends on the axis orientation only, $T_1 = T_3$ and $T_2 = T_4$; one can say that the hopping cannot distinguish forward steps from backward ones along a line. On the other hand, for $b-a$ odd the factor runs through all four values $1, \pm i, -1$, and all four $T_j$ are distinct, being related only through the hermiticity pairing. \\[.2em]
(b) To put it differently, it follows from \eqref{eq:Tj} that $\operatorname{tr} T_j = \sum_{a \in S} (-1)^a$ for all~$j$. Equal parities, $S = \{0,2\}$ and $\{1,3\}$,  give trace $\pm 2$, hence the matrices are, modulo the sign, twice \emph{orthogonal spin projectors}, $T_{1,3} = \pm(I + \sigma_1)$, $T_{2,4} = \pm(I - \sigma_1)$, which commute and are simultaneously diagonal, which suggests decoupling of horizontal and vertical motion. For opposite parities the trace vanishes, hence every $T_j$ is \emph{nilpotent}, $T_j^2 = 0$, the matrices $T_j$ are non-Hermitean and none two of them commute; we will see that the dynamics is then irreducibly two-component.
}
\end{remarks}
Using the matrices \eqref{eq:Tj} we can introduce the operator $\Heff^S$ on $\ell^2(\mathbb{Z}^2; \mathbb{C}^d)$ by
\begin{equation} \label{eq:Heff}
\bigl(\Heff^S \phi\bigr)(m,n) =
T_1\,\phi(m{-}1,n) + T_3\,\phi(m{+}1,n)
+ \ee^{-iBm}\,T_2\,\phi(m,n{+}1) + \ee^{iBm}\,T_4\,\phi(m,n{-}1)
\end{equation}
and to state the following duality claim:
\begin{theorem} \label{thm:duality}
Let $U_S$ be one of matrices \eqref{Kirch,antiKirch}--\eqref{pair_decoupl} and $B=2\pi\frac{p}{q}$ with coprime $p,q\in\N$, then for every $k \in (0,\infty) \setminus \pi\Z$ we have $k^2 \in \sigma(-\Delta_{B,U_S})$ if and only if $4\cos k \in \sigma\bigl(\Heff^S\bigr)$.
\end{theorem}
\begin{proof}
We will employ the appropriate Dirichlet-to-Neumann map. Identifying a lattice edge with the interval $[0,1]$ we suppose that $\psi$ solves $-\DD^2\psi=k^2\psi$ with constant (tangential component of the) potential $A$ there. The function is  uniquely determined by its endpoint values $\psi_0, \psi_1$, and
the outward quasi-derivatives at the endpoints are
\begin{equation} \label{eq:DtN}
\hspace{-1.2em} (\mathcal{D}_{\mathrm{out}}\psi)(0) = \frac{k}{\sin k}\bigl(\ee^{-iA}\psi_1 - \psi_0\cos k\bigr), \;
(\mathcal{D}_{\mathrm{out}}\psi)(1) = \frac{k}{\sin k}\bigl(\ee^{iA}\psi_0 - \psi_1\cos k\bigr).
\end{equation}
Indeed, the solution can be written as $\psi(x) = \ee^{iAx}(c\,\ee^{ikx} + d\,\ee^{-ikx})$ in which case $\mathcal{D}\psi = ik\,\ee^{iAx}(c\,\ee^{ikx} - d\,\ee^{-ikx})$. Identifying the endpoint values, we get $c + d = \psi_0$ and $c\,\ee^{ik} + d\,\ee^{-ik} = \ee^{-iA}\psi_1$; this system is solvable having determinant $-2i\sin k$, nonzero by assumption. This yields the uniqueness of the correspondence, and using the obtained coefficients, we get \eqref{eq:DtN}.

We will also employ the expansion bearing the name of Floquet, Bloch, or Gelfand, depending on the community to which you belong. Since by a simple gauge transformation one can check that integer multiples of the flux quantum do not affect the spectrum, both operators commute with translations on one cell vertically and $q$ cells horizontally, cf. \cite{BEL23} for the picture of the elementary cell with $q$ vertices and $4q$ `half-edges'. Consequently, one has to compare the corresponding fibers, differential one on the lattice graphs cell and the difference one on $\C^{dq}$, both having a discrete spectrum, and to show that
$$ 
\dim\ker\bigl(-\Delta_{B,U_S}(\theta) - k^2\bigr) = \dim\ker\bigl(\Heff^S(\theta) - 4\cos k\bigr)
$$ 
holds for any quasimomentum $\theta = (\theta_1, \theta_2) \in [0,2\pi)^2$, see, e.g., \cite[Sec.~4.2]{BK13}.

Let us begin with the map from the graph to the discrete lattice and consider a function $\psi \in \ker(-\Delta_{B,U_S}(\theta) - k^2)$. By \eqref{eq:projform}, the boundary values vector $\psi(v)$ at any vertex belongs to $E_+\C^4$; we associate with them the vector $\phi$ with the components  $\phi_a(v) := \langle \chi_a, \psi(v)\rangle$. The map $\psi \mapsto \phi$ is obviously linear and injective: $\phi=0$ implies $\psi=0$ by the Dirichlet-to-Neumann uniqueness mentioned above. We have to show that $\phi \in \ker(\Heff^S(\theta) - 4\cos k)$. Let us fix a vertex $v = (m,n)$ and index $a \in S$. As a subset of an orthonormal basis in $\C^4$, $\{\chi_b\}_{b\in S}$ is orthonormal, and $\sum_j \overline{(\chi_a)_j}\,\psi_j(v) = \phi_a(v)$ holds by the definition of $\phi_a(v)$, hence the identity $\langle \chi_a, \mathcal{D}\Psi(v)\rangle = 0$ can be using the first relation in \eqref{eq:DtN} rewritten as
\begin{equation} \label{eq:precond}
\sum_{j=1}^4 \overline{(\chi_a)_j}\,\omega_j(m)\,\psi(v^{\mathrm{next}}_j) = \cos k\:\phi_a(v),
\quad
\omega_1 = \omega_3 = 1, \;\; \omega_2(m) = \overline{\omega_4(m)} = \ee^{-iBm},
\end{equation}
where $v^{\mathrm{next}}_j$ is the neighboring vertex at the $j$th edge rooted in $v$ and where $\omega_j(m)$ are the Peierls phases \cite{Pe33}, written as $\ee^{\mp iA}$ in \eqref{eq:DtN}, coming from our choice of the Landau gauge. We observe that from the viewpoint of the vertex $v^{\mathrm{next}}_j$, the edge has to opposite orientation, hence $\psi(v^{\mathrm{next}}_j)$ equals $\sum_{b\in S}\phi_b(v_j)\,(\chi_b)_{\bar\jmath}\,$. Using then the identity \eqref{eq:reversal} together with $\overline{(\chi_a)_j}(\chi_b)_j = \tfrac14\, i^{(b-a)(j-1)}$, the left-hand side of \eqref{eq:precond} becomes $\tfrac14\sum_j \omega_j(m) \sum_b (T_j)_{ab}\, \phi_b(v_j)$, and multiplying the obtained relation by $4$ for a fixed $a\in S$, we arrive at $\bigl(\Heff^S(\theta)\phi\bigr)(v) = 4\cos k\,\phi(v)$, because the factors $\omega_j(m)$ are nothing but the phases appearing in \eqref{eq:Heff}.

In the opposite direction, consider a $\phi \in \ker(\Heff^S(\theta) - 4\cos k)$ and define $\psi(v) := \sum_{a\in S}\phi_a(v)\chi_a$ at every vertex; since $\sin k\ne 0$ by assumption, one can associate with these endpoint values the unique $k^2$-eigen\-function $\psi$ on each given edge. By construction, it satisfies $E_-\Psi(v) = 0$, and at the same time we have $E_+\mathcal{D}\Psi(v) = 0$, because using the above argument in the reverse order step by step, we obtain the corresponding discrete equation at $v$. Thus $\psi \in \ker(-\Delta_{B,U_S}(\theta) - k^2)$, so that the two maps are mutually inverse which concludes the proof.
\end{proof}

\subsection{Symmetries} \label{ss:symmetries}

An immediate consequence of Theorem~\ref{thm:duality} is that the spectrum of $-\Delta_{B,U_S}$ is $2\pi$-periodic in the momentum variable $k$. In fact, the spectrum is $\pi$-periodic.
\begin{proposition}\label{prop:chiral}
$\sigma(\Heff^S) = -\sigma(\Heff^S)$ holds for every $S$ and every flux, and consequently, $\sigma(-\Delta_{B,U_S})$ is $\pi$-periodic in $k$ on $(0,\infty)\setminus\pi\mathbb{Z}$.
\end{proposition}
\begin{proof}
For any fixed quasimomentum $\theta=(\theta_1,\theta_2)$, the vertical terms in the fiber $\Heff^S(\theta)$ become $\ee^{\mp i(Bm - \theta_2)}T_{2,4}$ flipping thus sign under the shift $\theta_2 \mapsto \theta_2 + \pi$, while the horizontal ones acquire the sign factor, $\phi_m \mapsto (-1)^m\phi_m$, between the cell boundaries so that they remain the same for even $q$ and require the additional shift $\theta_1 \mapsto \theta_1 + \pi$ for odd $q$. This means that $\Heff^S(\theta_1 + \pi q, \theta_2 + \pi)$ is unitarily equivalent to $-\Heff^S(\theta_1,\theta_2)$, therefore the spectrum of $\Heff^S$ is symmetric with respect to zero. Finally, $\cos(k+\pi) = -\cos k$, which implies the $\pi$-periodicity of all spectral bands of the original operator.
\end{proof}

Let us next look into further symmetries which would allow us to reduce the spectral analysis for the family \eqref{Kirch,antiKirch}--\eqref{pair_decoupl} to a smaller number of cases.
\begin{proposition}\label{prop:sign}
$\sigma(-\Delta_{B,-U_S}) = \sigma(-\Delta_{B,U_S})$ holds for every $S$ and every flux.
\end{proposition}
\begin{proof}
The spectrum of $-\Delta_{B,U_S}(\theta)$ away from $\pi\mathbb{Z}$ can be determined using the vertex scattering matrix which is in the scale-invariant case independent of energy, $S(k)=U_S$. The eigenvalues then come from solutions of the secular equation \cite[Sec.~2.2]{BK13}, \cite[Sec.~16.2]{Ku24}, which reads $\det(I-\ee^{ikL}\,\mathbb{S}\,\PP)=0$, where $\mathbb{S}(k)$ is the bond scattering matrix, $L$ has edge lengths on the diagonal and $\PP$ carries the edge reversal (connecting two bonds of each edge) together with the Peierls and Floquet phase factors. In our case all edges are of length one, and replacing $U_S$ by $-U_S$ one changes $\mathbb{S}$ to $-\mathbb{S}$ and $\det(I-\ee^{ikL}\,\mathbb{S}\,\PP)$ to $\det(I+\ee^{ik}\mathbb{S}\mathcal{P}) = \det(I-\ee^{i(k+\pi)}\mathbb{S} \mathcal{P})$, thus the spectrum of $-\Delta_{B,-U_S}(\theta)$ is shifted by $\pi$ in $k$, however, by Proposition~\ref{prop:chiral} it is $\pi$-periodic, hence unchanged, and the conclusion carries over to the direct integral over $\theta$. The exceptional set $\pi\mathbb{Z}$ is invariant under the shift.
\end{proof}

This is not the only symmetry, another one comes from complex conjugation which is, in view of the unitarity, equivalent to transposition.
\begin{proposition}\label{prop:conj}
$\sigma(\Delta_{B,\bar{U}_S}) = \sigma(-\Delta_{B,U_S})$ holds for every $S$ and every flux.
\end{proposition}
\begin{proof}
Complex conjugation is an antiunitary equivalence between $-\Delta_{B,U}$ with the plaquette flux $\Phi$ and $-\Delta_{-B,\bar U}$ with the flux $-\Phi$ and we have $\bar U_S = U_{-S\,(\!\bmod\,4)}$. We need to care about nonreal $U_S$ only, that is, the pairs \eqref{single1} and \eqref{single3} for $d\ne 2$, and  \eqref{double_side} and \eqref{double_middle} for $d=2$. Each of these two groups of four is invariant with respect to the swap $S\leftrightarrow -S\,(\bmod\,4)$ and within each we have symmetry with respect to the flux reversal, $p\leftrightarrow q-p$ (recall that integer multiples of the flux quantum do not matter).
\end{proof}

\section{The main result} \label{s:main}
\setcounter{equation}{0}

The described symmetries reduced the number of different situations from fourteen to five. It appears that from the spectral point of view the true number of isospectral classes is three:
 \begin{enumerate}[I.]
 \setlength{\itemsep}{3pt}
 \item Eight couplings \eqref{Kirch,antiKirch}--\eqref{single3} with $d=1$ and $d=3$.
 \item Four couplings \eqref{double_side} and \eqref{double_middle} with $d=2$.
 \item The pair of couplings \eqref{pair_decoupl} with $d=2$.
 \end{enumerate}
\begin{theorem}\label{thm:class}
The spectrum of $-\Delta_{B,U_S}$ for the coupling matrices $U_S$ of the indicated classes looks as follows:
\begin{enumerate}[I.]
 \setcounter{enumi}{0}
 \setlength{\itemsep}{3pt}
 \item For every rational value of the flux, $\Phi= 2\pi\frac{p}{q}$, it coincides with that of the Kirchhoff lattice, i.e.  $k^2 \in \sigma(-\Delta_{B,U_S})$ if and only if $4\cos k \in \sigma(H_{p/q})$, where $H_{p/q}$ is the Harper operator, together with eigenvalues of infinite multiplicity (flat bands) at every $k \in \pi\mathbb{Z}$.
 \item For every rational flux, $\Phi= 2\pi\frac{p}{q}$, the spectrum is by Theorem~\ref{thm:duality} dual that of the two-component operator with hopping matrices $T_{1,3} = \sigma_3 \mp i\sigma_2$ and  $T_{2,4} = \sigma_3 \pm i\sigma_1$ (for the given choice of the basis in $E_+\C^4$), or equivalently to  that of the dimerized Harper chain
\begin{equation} \label{eq:chain}
\begin{aligned}
b_{m-1} + \sin\xi_m\,a_m + \cos\xi_m\,b_m & = 2\cos k\;a_m, \\[.1em]
a_{m+1} + \cos\xi_m\,a_m - \sin\xi_m\,b_m & = 2\cos k\;b_m,
\end{aligned}
\end{equation}
where $\xi_m = Bm - \theta_2$. There are no flat bands; the butterfly has three fractal clusters per period, centered at $k \equiv \tfrac\pi4, \tfrac\pi2, \tfrac{3\pi}4 \pmod\pi$.
 \item The lattice graph decouples into independent horizontal and vertical one-dimensional lines, gauge equivalent to the non-magnetic situation; the spectrum is $[0,\infty)$ independently of the flux, and there are no flat bands.
\end{enumerate}
\end{theorem}
\begin{proof}
\emph{Class I.} For $d = 1$ the hopping matrices are just numbers; if $S = \{m\}$, formula \eqref{eq:Tj} gives $T_j = (-1)^m$ for all four directions, hence the dual operator is the $(-1)^m$ multiple of the well-known two-dimensional Harper operator $H_\alpha,\,\alpha=\frac{p}{q}$, cf.~\cite{JKK22}, also known under other names such as the discrete magnetic Laplacian \cite{Sh94}. In view of the chiral symmetry from Proposition~\ref{prop:chiral}, $\sigma(\Heff^S) = \sigma(-\Heff^S)$, so that the spectrum is the same for both signs, identical with that of the Kirchhoff vertex coupling to which the first of the matrices \eqref{Kirch,antiKirch} corresponds. For the other half of this class, $d=3$, the coupling is of the form $-U_{\{m\}}$ so the result follows from Proposition~\ref{prop:sign}. As for the flat bands, in the Kirchhoff case they are the usual Dirichlet loop states having for rational flux zeros at the vertices. For the other cases one can check their existence numerically, see Remark~\ref{rem:flat} below.

\emph{Class II.} For $S = \{0,3\}$ formula \eqref{eq:Tj} gives $T_3 = \sigma_3 + i\sigma_2$, $T_2 = \sigma_3 + i\sigma_1$ and their adjoints. In the Landau gauge we use, $\Heff$ is diagonalized by Bloch waves in the $n$-direction: writing $\phi(m,n) = \ee^{i\theta_2 n}\varphi_m$ and rotating the spinor by the Hadamard matrix $W = (\sigma_1 + \sigma_3)/\sqrt{2}$, the components $\tilde\varphi_m = (a_m, b_m)^{\mathsf T} = W \varphi_m$ satisfy the chain equation \eqref{eq:chain}. The other three members are obtained by conjugation and sign flip, hence they are isospectral by Propositions~\ref{prop:sign} and \ref{prop:conj}. The absence of flat bands
and the shape of the spectrum in dependence on the flux shown in Figure~\ref{newbutterfly} for $q=1,\dots,30$ have been obtained numerically. The clusters centers at $\lambda = 0, \pm2\sqrt2$ are the loci of the band contacts of the flux-one-half line.

\emph{Class III.} Consider $S=\{0,2\}$ in which case $E_+\C^4 = \operatorname{span}\{\chi_0, \chi_2\} = \{\psi : \psi_l = \psi_r,\ \psi_u = \psi_d\}$; these coincidences of the boundary values mean that the wave function is continuous along each straight line through the vertex, with no relation between the two lines.
The derivative conditions $\langle\chi_0, \DD\Psi\rangle = \langle\chi_2, \DD\Psi\rangle = 0$ from \eqref{eq:projform} are equivalent to $\DD_l + \DD_r = 0$ and $\DD_u + \DD_d = 0$, i.e. to Kirchhoff conditions separately on each line. In the operator sense, the lattice is thus a disjoint union of free one-dimensional `wires'; the tangential component vector potential on a straight wire is easily removed by a gauge transformation, so each wire contributes by $[0,\infty)$ which is also the total spectrum, independent of the flux. Equivalently, on the dual side relation \eqref{eq:Tj} gives $T_{1,3} = I + \sigma_1$, $T_{2,4} = I - \sigma_1$, which commute; in the $\sigma_1$ eigenbasis the fibers decouple into the free horizontal chain with eigenvalues $4\cos\frac{\theta_1 + 2\pi l}{q}$ and the on-site vertical values $4\cos(Bm - \theta_2)$, whose union over $\theta$ is $[-4,4]$. The other case, $S=\{1,3\}$, is covered by Proposition~\ref{prop:sign}.
\end{proof}
\begin{remarks} \label{rem:flat}
{\rm Except Kirchhoff, the existence of flat bands for class I has been checked numerically, see the disclosure claim below. Remarkably, for $d=3$ there is a way to see their origin. Take. e.g., $S = \{0,2,3\}$, where formula \eqref{eq:Tj} gives
$$
T_1 = {\scriptsize \begin{pmatrix} 1 & 1 & -1\\[1pt] 1 & 1 & -1\\[1pt] 1 & 1 & -1\end{pmatrix}},
\quad
T_2 = {\scriptsize \begin{pmatrix} 1 & -1 & i\\[1pt] -1 & 1 & -i\\[1pt] i & -i & -1\end{pmatrix}},
\quad
T_3 = T_1^*,\quad T_4 = T_2^* = \overline{T_2},
$$
in the basis $(\chi_0, \chi_2, \chi_3)$; the other three $d = 3$ patterns are obtained analogously. By the symmetries, the above proof shows that spectrally these cases coincide with Harper's operator, but on the level of wave functions they are different. The structure of hopping matrices shows that all four directions are now distinct; their traces are $\sum_{a\in S}(-1)^a = -(-1)^{m'} = \pm1$, where $m'$ is the index missing in $S$ and they satisfy $T_j^2 = \pm T_j$ being neither nilpotent as in class II nor Hermitean as in class III. Concerning the flat bands, each fiber component spectrum consists of $3q$ eigenvalues of which $q$ coincide with those of $-\sigma(H_{p/q}(\theta))$ and the two others, each of multiplicity $q$, are situated at $\pm 4$, the edges of the duality window. \\[.2em]
(b) The collapse from five to three spectral classes has two sources which the symmetries alone do not reveal. The first is the exact $\pi$-periodicity from Proposition~\ref{prop:chiral}, which converts the $k \mapsto k+\pi$ shift of the sign flip into an equality of spectra, merging thereby each $d=3$ coupling with its $d=1$ partner; the second is the bipartite gauge equivalence of the four scalar duals, which merges the three class I symmetry cases.}
\end{remarks}

\section{Discussion} \label{s:discussion}
\setcounter{equation}{0}

\subsection{The present result} \label{ss:present}

\begin{figure}[htb]
\centering
\includegraphics[width=\textwidth]{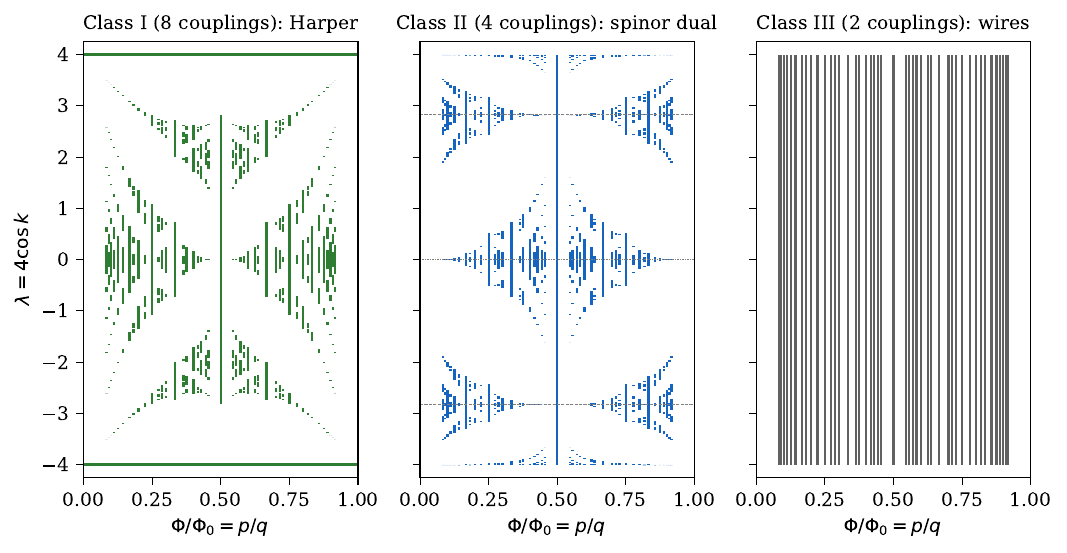}
\caption{The three spectral classes of scale-invariant circulant couplings, in the duality variable $\lambda = 4\cos k$, for all coprime $p/q$ with $q \le 12$. In class I the heavy lines at $\lambda = \pm4$ indicate the flat bands at $k \in \pi\mathbb{Z}$, in class II the dashed lines at $\lambda = \pm2\sqrt2$ and $\lambda = 0$ mark the three fractal clusters centers.}
\label{fig:classes}
\end{figure}

The three spectral classes are shown in Figure~\ref{fig:classes} in terms of the duality variable $\lambda = 4\cos k$; in view of Theorem~\ref{thm:duality} and Proposition~\ref{prop:chiral}, each period $[n\pi, (n+1)\pi]$ of the momentum axis carries one copy of the respective diagram. Since the aim here is a comparison, the resolution is coarser than in Figure~\ref{newbutterfly}.
The differences between the two nontrivial classes are obvious, in the periodicity, shape of the fractal clusters, or the fact that class II is gapless at the half-of-the-quantum flux. To elucidate the mechanism behind this difference, it instructive to inspect the dispersion in the zero-field situation.
\begin{corollary}\label{cor:disp}
For $B = 0$ the class-II lattice has the dispersion relation
\begin{equation}
\cos 2k \;=\; \cos\theta_1\,\cos\theta_2 .
\label{eq:disp}
\end{equation}
In particular, the spectrum is all of $[0,\infty)$, with no flat bands, and the two dispersion surfaces of $-\Delta_{0,U_S}$ touch \emph{conically} at the quasimomenta $(\theta_1,\theta_2) = (0,\pi), (\pi,0)$, that is, at the momenta $k \equiv \tfrac\pi2 \pmod\pi$.
\end{corollary}
\begin{proof}
For definiteness, consider $S = \{0,3\}$, the other cases are similar. If $B = 0$, the operator \eqref{eq:Heff} is invariant under all lattice translations, so the Fourier Ansatz $\phi(m,n) = e^{i(\theta_1 m + \theta_2 n)}\varphi$ gives the fiber a $2\times2$ matrix representation, namely
\begin{align*}
H(\theta) &= (\sigma_3 - i\sigma_2)\,\ee^{-i\theta_1} + (\sigma_3 + i\sigma_2)\,\ee^{i\theta_1}
+ (\sigma_3 + i\sigma_1)\,\ee^{i\theta_2} + (\sigma_3 - i\sigma_1)\,\ee^{-i\theta_2}\\[2pt]
&= 2(\cos\theta_1 + \cos\theta_2)\,\sigma_3 \;-\; 2\sin\theta_1\,\sigma_2 \;-\; 2\sin\theta_2\,\sigma_1
\:=:\: \mathbf h(\theta)\cdot\boldsymbol\sigma .
\end{align*}
The matrix is traceless -- recall that $\operatorname{tr} T_j = 0$ for class II by Remark~\ref{hop_differences}(b) -- hence its eigenvalues are $\pm|\mathbf h(\theta)|$, for which a simple manipulation gives
$$
|\mathbf{h}(\theta)|^2 = 4(\cos\theta_1 + \cos\theta_2)^2 + 4\sin^2\theta_1 + 4\sin^2\theta_2
= 8\,\bigl(1 + \cos\theta_1\cos\theta_2\bigr).
$$
The spectral duality from Theorem~\ref{thm:duality} gives $16\cos^2 k = 8(1 + \cos\theta_1\cos\theta_2)$, i.e.\ $2\cos^2 k - 1 = \cos\theta_1\cos\theta_2$, which is \eqref{eq:disp}. The right-hand side range is $[-1,1]$ as $\theta$ runs over the Brillouin zone, so every $k > 0$ is attained and the spectrum is $[0,\infty)$. The dispersion surfaces touch at the points where $\mathbf{h}(\theta) = 0$: expanding around $(\theta_1, \theta_2) = (0,\pi)$ with $\theta_1 = \delta_1$, $\theta_2 = \pi + \delta_2$, one finds $|\mathbf{h}(\theta)| = 2\sqrt{\delta_1^2 + \delta_2^2} + O(\delta^3)$, a conical (Dirac) contact at $\lambda = 4\cos k = 0$, i.e.\ at $k \equiv \tfrac\pi2 \pmod\pi$.
\end{proof}

The dispersion \eqref{eq:disp} contrasts with $4\cos k = 2\cos\theta_1 + 2\cos\theta_2$, its Kirchhoff counterpart, which contains the sum of the quasimomentum cosines instead of their product; the two spectra are gapless but the Kirchhoff's contains flat bands which are absent here. The conical contacts of Corollary~\ref{cor:disp} can be regarded as the seeds of the central cluster of the class-II butterfly; one can recall the central cluster of the graphene butterfly emanating from the Dirac points of the honeycomb dispersion. Here the cluster centered at $\lambda = 0$ in Figure~\ref{fig:classes} grows out of the two zero-field Dirac contacts at $k \equiv \tfrac\pi2$.
\begin{remark}\label{rem:squared}
{\rm To indicate to what degree does the mechanism suggested by Corollary~\ref{cor:disp} extends to a nonzero flux, consider the square of the chain operator $H_c$ acting on $\psi={a_m\choose b_m}$ according to \eqref{eq:chain}. A direct computation gives
\begin{equation} \label{eq:squared}
\begin{aligned}
(H_c^2\psi)_{a_m} &= 2\,a_m + \cos\xi_m\, a_{m+1} + \cos\xi_{m-1}\, a_{m-1} + \rho_m\, b_{m-1},\\
(H_c^2\psi)_{b_m} &= 2\,b_m + \cos\xi_{m+1}\, b_{m+1} + \cos\xi_m\, b_{m-1} + \rho_{m+1}\, a_{m+1},
\end{aligned}
\end{equation}
where $\rho_m := \sin\xi_m - \sin\xi_{m-1}$, because the mixed same-cell terms cancel. We note that the diagonal is independent of $m$; putting thus $\epsilon=2\cos k$ and using the identity $\epsilon^2-2=2\cos 2k$, we see that the operator $G := \tfrac12(H_c^2 - 2)$ satisfies the `doubled' duality relation: $k^2 \in \sigma(-\Delta_{B,U_S})$ if and only if $\cos 2k\in\sigma(G)$. Within each sublattice, $G$ is a purely off-diagonal Harper chain with no on-site term and hopping amplitudes $\tfrac12\cos\xi_m$ between consecutive $a$'s, and similarly for the $b$'s. The two sublattices are coupled by the bond $\rho_m$ which can be rewritten as
$$
\rho_m = 2\sin\!\Bigl(\frac{\pi p}{q}\Bigr)\,\cos\!\Bigl(\xi_m - \frac{\pi p}{q}\Bigr),
$$
being a Harper-modulated amplitude of the same frequency and the global strength $2\sin(\pi p/q)$. It vanishes identically if and only if $B = 0$, and in that case the sublattices decouple and the Bloch reduction of the $a$-chain reproduces dispersion relation \eqref{eq:disp}. At a nonzero flux, this term prevents such a scalar reduction and generates the central cluster of the butterfly. Note that the amplitude is maximal at half-of-the-quantum flux when the spectrum is gapless. }
\end{remark}

On another topic, as the motivation to this paper we mentioned the spectral problem for the same lattice with a simple vertex coupling that violated the time-reversal invariance \cite{BEL23}. The said coupling was not scale-invariant so the obtained spectrum lacked the $k$-periodicity; it was distorted at low energies while in the high-energy asymptotics it approached the Hofstadter butterfly pattern, except for narrow regions around the values where the Kirchhoff lattice had the flat bands. Here the result is significantly different, and not only because of the periodicity, natural in our case; the main, and to a degree surprising distinction is the way in which the time asymmetry, in the present model always maximal by Remark~\ref{maxviol}, enters. The new spectrum of class II is related to  the time-reversal invariance violation as none of the matrices \eqref{double_side} and \eqref{double_middle} is transpose invariant. On the other hand, the spectrum of class I has nothing to do with the time-reversal invariance, or the absence of it; one half of the couplings \eqref{Kirch,antiKirch}--\eqref{single3} gives rise to time-reversal invariant dynamics, the other to an irreversal one.

The result opens various questions. The most natural one concerns the nature of the spectrum for irrational values of the flux; one can conjecture that it would be a Cantor set in analogy with the answer to the classical Ten Martini Problem~\cite{AJ09}. Of interest is also whether the Haussdorff dimension dichotomy known in the critical almost Mathieu situation \cite{BLQZ19, LS16} carries over to the class II spectrum. At the same time, the present result certainly does not exhaust all questions about the rational flux case; one can ask, for instance, what would be a counterpart of the Thouless conjecture \cite{BK95, JKK22, Th83, Th90} here.

\subsection{Related problems in the literature} \label{ss:related}

Our construction is connected with several areas beyond the familiar Kirchhoff and almost-Mathieu correspondence. \\[.2em]
\emph{Dualities:} The differential-difference operator connection is the core element of our construction. The duality idea can be traced back to \cite{vB85, Ca97} for weighted Kirchhoff on equilateral graphs and to \cite{Ex97} for $\delta$- and $\delta'$-couplings; an abstract boundary-triple formulation is due to \cite{Pa06, Pa12}. Theorem~\ref{thm:duality} extends this to scale-invariant couplings, for which the dual operator is genuinely matrix-valued, of fiber dimension $\dim\ker(U - I)$. In this connection one can mention results on honeycomb lattices. Kuchment and Post \cite{KP07} reduced the non-magnetic graphene quantum graph to a discrete two-component problem, Becker, Han, and Jitomirskaya \cite{BHJ19} carried out the magnetic reduction and proved Cantor spectrum for irrational flux, see also \cite{BZ19}. The difference is that the honeycomb is an overlay of two lattices; in our case the matrix character comes from properties of a single vertex. \\[.1em]
\emph{Discrete magnetic operators on periodic graphs:} The operators \eqref{eq:Heff} can be naturally regarded as magnetic Schr\"odinger operators on $\mathbb{Z}^2$ with a $d$-dimensional fiber, $d=1,2,3$. A general band structure analysis for $\Z^d$ periodic discrete graphs was done by Korotyaev and Saburova, both in the non-magnetic \cite{KS14} and magnetic \cite{KS17} cases. From the similar reason as above, however, it does not apply directly to our present problem because the multidimensionality is not there a single-site matter but comes from the number of vertices in the elementary cell.  \\[.1em]
\emph{Harper-type models with modulated hopping:} One-dimensional quasiperiodic operators in which the hopping exhibits a (typically incommensurate) modulation appear repeatedly in the literature; one can recall, for instance, the extended Harper model expressed in terms of off-diagonal cosine factors analyzed  by Avila, Jitomirskaya, and Marx \cite{AJM17}. For comparison, squaring the chain \eqref{eq:chain} and projecting onto one sublattice produces, by \eqref{eq:squared}, up to the coupling $\rho_m$ of strength $2\sin(\pi p/q)$ made explicit in Remark~\ref{rem:squared}, a purely off-diagonal Harper chain with hopping $\cos\xi_m$ and constant diagonal. However, the chain itself does not belong to the exactly solved extended-Harper family; the modulation appears both in the dimerized bond and the staggered site with the unitary constraint $\cos^2\xi_m + \sin^2\xi_m = 1$,  which recalls rather a structure appearing in the literature on the unitary side --- in CMV matrices and the unitary almost-Mathieu operator \cite{FOZ17}, and in discrete-time quantum walks in electric and magnetic gauge fields \cite{CGWW19, CRW13}. In contrast, here we deal with a Hermitean lattice operator having a unitary parametrization.  \\[.1em]
\emph{Hofstadter-type spectra.} Butterfly patterns coming from multicomponent lattice models appear in various places; one can mention. e.g.,  Rammal's spectrum of the honeycomb (graphene) tight-binding model in a field \cite{AEG14, KS96, Ram85}, butterflies of flat-band lattices \cite{AAM96},  deformed square-lattice butterflies with next-nearest-neighbour hopping \cite{HK90}, or Hofstadter problems with non-Abelian gauge fields \cite{OBSZL05}. Our class II is distinct from all of these in that the magnetic field is Abelian and the matrix structure of the hopping is not a gauge degree of freedom but a fixed, rank-one nilpotent imprint of the vertex coupling --- each bond transmits a single scalar channel, the hallmark of the quantum-graph origin.

\subsection*{AI usage disclosure}

The text was prepared in dialogue with Anthropic Claude Fable used to check some steps in the reasoning, perform the numerical computations, and assist in literature search. The AI suggestions were verified and the text was written by the authors who take the full responsibility for its correctness.

\subsection*{Data availability statement}

Data are available in the article.

\subsection*{Conflict of interest}

The authors have no conflict of interest.


\end{document}